\documentclass[11pt]{article}

\usepackage{amssymb}
\usepackage[english]{babel}
\usepackage{graphicx}
\usepackage{latexsym}

\usepackage{amsmath,amsfonts,amssymb,latexsym,graphics,epsfig,url}
\usepackage{color}
\usepackage{amsthm}
\usepackage[english]{babel}
\usepackage{graphicx}
\usepackage{soul}
\usepackage[utf8]{inputenc}

\newtheorem{theorem}{Theorem}[section]

\def\Cay{\mbox{\rm Cay}}

\DeclareMathOperator{\Irep}{Irep}
\DeclareMathOperator{\spec}{sp}

\def\Z{\ns Z}

\def\j{\mbox{\boldmath $j$}}

\def\z{\mbox{\boldmath $z$}}

\def\vecalpha{\mbox{\boldmath $\alpha$}}
\def\vecbeta{\mbox{\boldmath $\beta$}}
\def\vecgamma{\mbox{\boldmath $\gamma$}}
\def\vecdelta{\mbox{\boldmath $\delta$}}

\def\vecx{\mbox{\boldmath $x$}}

\def\vecphi{\mbox{\boldmath $\phi$}}
\def\vec0{\mbox{\boldmath $0$}}

\def\A{\mbox{\boldmath $A$}}
\def\B{\mbox{\boldmath $B$}}

\def\D{\mbox{\boldmath $D$}}

\def\G{\Gamma}
\def\I{\mbox{\boldmath $I$}}
\def\J{\mbox{\boldmath $J$}}

\def\U{\mbox{\boldmath $U$}}

\def\Z{\ns{Z}}
\def\I{\mbox{\boldmath $I$}}
\def\J{\mbox{\boldmath $J$}}

\def\1{\mbox{\boldmath $1$}}

\def\CC{\mathbb C}

\def\G{\Gamma}

\def\Z{\mathbb Z}

\begin{document}
	\title{Token graphs of Cayley graphs as lifts
		\thanks{The research of the first two authors has been supported by
			AGAUR from the Catalan Government under project 2021SGR00434 and MICINN from the Spanish Government under project PID2020-115442RB-I00.
			The second author's research was also supported by a grant from the  Universitat Polit\`ecnica de Catalunya with references AGRUPS-2022 and AGRUPS-2023. The third and the fourth authors gratefully acknowledge support for this
			research from APVV Research Grants 19-0308 and 22-0005, and VEGA Research
			Grants 1/0567/22 and 1/0069/23.}}

\author{C. Dalf\'o$^a$, M. A. Fiol$^b$,  S. Pavl\'ikov\'a$^c$, and J. \v{S}ir\'a\v{n}$^d$\\
	\\
	{\small $^a$Departament. de Matem\`atica}\\ {\small  Universitat de Lleida, Lleida/Igualada, Catalonia}\\
	{\small {\tt cristina.dalfo@udl.cat}}\\
	{\small $^{b}$Departament de Matem\`atiques}\\ 
	{\small Universitat Polit\`ecnica de Catalunya, Barcelona, Catalonia} \\
	{\small Barcelona Graduate School of Mathematics} \\
	{\small  Institut de Matem\`atiques de la UPC-BarcelonaTech (IMTech)}\\
	{\small {\tt miguel.angel.fiol@upc.edu}}\\
	{\small $^{c}$ Inst. of Information Engineering, Automation, and Math, FCFT}, \\
	{\small Slovak Technical University, Bratislava, Slovakia}\\
	{\small {\tt sona.pavlikova@stuba.sk} }\\
	{\small $^{d}$Department of Mathematics and Statistics}\\
	{\small The Open University, Milton Keynes, UK}\\
	{\small Department of Mathematics and Descriptive Geometry}\\
	{\small Slovak University of Technology, Bratislava, Slovakia}\\
	{\small {\tt jozef.siran@stuba.sk}}
}

	\maketitle
	
	\begin{abstract}
		This paper describes a general method for representing $k$-token graphs of Cayley graphs as lifts of voltage graphs. This allows us to construct line graphs of circulant graphs and Johnson graphs as lift graphs on cyclic groups. As an application of the method, we derive the spectra of the considered token graphs.
		The method can also be applied to dealing with other matrices, such as the Laplacian or the signless Laplacian,
		and to construct token digraphs of Cayley digraphs.
	\end{abstract}

	\section{Introduction}
	
	A considerable amount of research in algebraic graph theory has been devoted
	to constructions of larger graphs from smaller ones, enabling the determination of
	properties of larger graphs in an algebraically controllable way from
	properties of the `parent' graphs. Although examples of such constructions
	abound,
	we will focus only on two of them, namely, the token graph construction and
	the lifting construction. They turn out to be related in an interesting way,
	and before an explanation, we offer an informal description of the two
	constructions.
	
	Given a graph $G=(V,E)$ on $n$ vertices, its $k$-token graph $F_k(G)$ has a vertex set consisting of all the ${n\choose k}$ configurations of $k$ indistinguishable tokens placed in different vertices of $G$, and two vertices are adjacent if one is obtained from the other by moving one token along an edge of $G$. These graphs were also called {\em symmetric $k$-th power of a graph} by Audenaert, Godsil, Royle, and Rudolph \cite{agrr07}, and {\em $k$-tuple vertex graphs} in  Alavi, Behzad, Erd\H{o}s, and Lick \cite{abel91}.
	Token graphs have applications in physics (quantum mechanics) and in the graph isomorphism problem (because invariants of the $k$-token $F_k(G)$ are also invariants of $G$). For more details, see \cite{agrr07}.
	Some properties of token graphs were studied by Fabila-Monroy, Flores-Pe\~{n}aloza, Huemer,  Hurtado,  Urrutia, and  Wood \cite{ffhhuw12}, including the connectivity, diameter, clique and chromatic numbers, and Hamiltonian paths.
	In particular, the connectivity of token graphs of trees was studied by Fabila-Monroy, Leaños,  and Trujillo-Negrete
	in \cite{flt22}.
	
	For an informal description of the lifting construction, think of a graph $G$
	(a {\em base graph}) endowed by an assignment $\alpha$ of elements of a group
	$\Gamma$ on arcs of $G$ (a {\em voltage assignment}). The pair $(G,\alpha)$
	gives rise to a {\em lift} $G^\alpha$, a larger graph that may be thought of
	as obtained by `compounding' $|\Gamma|$ copies of $G$, joined in-between in a
	way dictated by the voltage assignment. A general necessary and sufficient
	condition for a graph $H$ to arise as a lift of a (smaller) graph is the existence
	of a subgroup of the automorphism group of $H$ with a free action on vertices
	of $H$; see Gross and Tucker \cite{gt87}. Lifts have found numerous applications
	in areas of graph theory that are as versatile as the degree/diameter problem on the
	one hand and the Map Color Theorem on the other hand.
	In the first case, the diameter of the lift can be conveniently expressed
	in terms of voltages on the edges of the base graph. Besides its theoretical importance,
	this fact can be used to design efficient diameter-checking algorithms, see Baskoro, Brankovi\'c, Miller, Plesn\'{\i}k,  Ryan, and  \v{S}ir\'a\v{n} \cite{bbmprs97}.
	In the context of the degree-diameter problem, we address the interested reader to the comprehensive survey by Miller and  \v{S}ir\'a\v{n} \cite{ms13}.
	Other prominent examples of
	applications also include the use of lift graphs in the study of the automorphisms of $G^{\alpha}$ and in Cayley graphs, which are lifts of one-vertex graphs
	(with loops and semi-edges attached).
	
	The surprising way the two constructions turn out to be related is based on
	the observation that {\sl automorphism groups of $k$-token graphs of Cayley
		graphs admit subgroups acting freely on vertices} and, hence, they can be described
	as lifts of smaller graphs by voltage assignments. This enables one to link
	the study of $k$-token graphs, Cayley graphs, and lifts in a novel way by
	introducing a general method to represent $k$-token graphs of Cayley graphs as
	lifts of voltage graphs. In particular, this results in representations of
	line graphs of circulant graphs and Johnson graphs as lift of smaller graphs,
	with voltages in cyclic groups. As another application of this method, we
	derive the spectra of the token graphs considered, and generate infinite families
	of graphs with given maximum degree and eigenvalues contained in a certain 
	interval.

	This paper is structured as follows. In Section \ref{sec:preliminaries}, there are the preliminaries, some notation, and some known results. In Section \ref{sec:Johnson}, we construct line graphs of circulant graphs and Johnson graphs as lift graphs on the cyclic group. Applying this method, we derive the spectra of the considered token graphs in Section \ref{sec:cay}.
	Finally, in Section \ref{sec:further}, we explain how to extend our results in the following directions: computing eigenspaces, working with the universal matrix instead of the adjacency matrix, and working with digraphs.
	
	\section{Preliminaries}
	\label{sec:preliminaries}
	
	Let $G$ be a graph with vertex set $V(G)=\{1,2,\ldots,n\}$ and edge set $E(G)$. 
	If necessary, we consider every edge $\{u,v\}$ as a `digon' formed by the arcs (or directed edges) $a^+=(u,v)$ and $a^-=(v,u)$. For a given integer $k$ such that $1\le k \le n$, the {\em $k$-token graph} $F_k(G)$ of $G$ is the graph in which the vertices of $F_k(G)$ correspond to $k$-subsets of $V(G)$. Two vertices are adjacent when the corresponding subsets' symmetric difference are the edge's end-vertices in $E(G)$.
	In particular, notice that
	$F_1(G)=G$ and, by symmetry,
	$F_k(G)=F_{n-k}(G)$. 
	
	Moreover,
	$F_k(K_n)$  is the Johnson graph $J(n,k)$. For example, Figure \ref{fig2} shows the Johnson graph $J(5,2)\cong F_2(K_5)$.
	In general, the Johnson graph $J(n,k)$, with $k\le n-k$, is a distance-regular graph with degree $k(n-k)$, diameter $d=k$, and for $j=0,1,\ldots,d$, intersection parameters
	$$
	b_j=(k-j)(n-k-j), \qquad c_j=j^2,
	$$
	eigenvalues
	$$
	\lambda_j=(k-j)(n-k-j)-j,
	$$
	and multiplicities
	$$
	m_j={n\choose j}-{n\choose j-1}.
	$$
	For more information about Johnson graphs, see Godsil \cite{g93}.
	
	In our study, some infinite families of line graphs appear. 
	Concerning their spectra, recall that if a regular graph $G$ with $n$ vertices and $m$ edges has the spectrum
	$$
	\spec(G)=\{k,\lambda_1^{[m_1]},\ldots,\lambda_d^{[m_d]}\},
	$$
	then its line graph $L(G)$ has spectrum
	$$
	\spec(L(G))=\{2k-2,(k-2+\lambda_1)^{[m_1]},\ldots,(k-2+\lambda_d)^{[m_d]},-2^{[m-n]}\}.
	$$
	The graphs with the least eigenvalue not smaller than $-2$ are completely identified in three categories:
	the line graphs of bipartite graphs,
	the generalized line graphs, and a  finite class of graphs arising from  root systems, see Biggs \cite{biggs}, and 
	Cameron,  Goethals, Seide, and  Shult \cite{cgss76}.
	
	Let $\G$ be a group. An ({\em ordinary\/}) {\em voltage assignment} on the graph $G=(V,E)$  is a mapping $\alpha: E\to \G$ with the property that $\alpha(a^-)=(\alpha(a^+))^{-1}$ for every pair of opposite arcs $a^+,a^-\in E$. Then, the graph $G$ and the voltage assignment $\alpha$ determine a new graph $G^{\alpha}$, called the {\em lift} of $G$, which is defined as follows. The vertex and arc sets of the lift are simply the Cartesian products $V^{\alpha}=V\times \G$ and $E^{\alpha}=E\times \G$, respectively. Moreover, for every arc $a\in E$ from a vertex $u$ to a vertex $v$ for any $u,v\in V$ (possibly, $u=v$) in $G$, and for every element $g\in \G$, there is an arc $(a,g)\in E^{\alpha}$ from the vertex $(u,g)\in V^{\alpha}$ to the vertex $(v,g\alpha(a))\in V^{\alpha}$.
	The interest of this construction is that, from the base graph $G$ and the voltages, we can deduce some properties of the lift graph $G^{\alpha}$. This is usually done through a matrix associated with $G$ that, in the case when the group $\G$ is Abelian of rank $r$, every entry of such a matrix is a polynomial in $r$ variables $z_1,\ldots,z_r$ with integer coefficients. So it is called 
	the associated `polynomial matrix' $\B(\z)$, where $\z=(z_1,\ldots,z_r)$. 
	If, in $G$, there are $r$ arcs from $u$ to $v$ with voltages $j_1,\ldots,j_r$, the entry $uv$ of $\B(z)$
	has a term of the form $z_1^{j_1}\cdots z_r^{j_r}$
	(see examples with $r=1$ and $r=2$ in the following sections).
	Then, the whole (adjacency or Laplacian) spectrum and eigenspaces of $G^{\alpha}$ can be retrieved from $\B(\z)$. 
	More precisely, Dalf\'o, Fiol, Miller, Ryan, and \v{S}ir\'a\v{n} \cite{dfmrs17} proved the following result.
	\begin{theorem}[\cite{dfmrs17}]
		\label{th:sp-lifts}
		Let $R(n)$ be the set of $n$-th roots of unity, and consider the base graph $G=(V,E)$ with voltage assignment $\alpha$ on the Abelian group  $\G=\Z_{n_1}\times\cdots\times \Z_{n_r}$. If $\vecx=(x_u)_{u\in V}$ is an eigenvector of $\B(\z)$ and $z_i\in R(n_i)$ with eigenvalue $\lambda$, then the vector $\vecphi=(\phi_{(u,\j)})_{(u,\j)\in V^{\alpha}}$ with 
		$\j=(j_1,\ldots,j_r)\in \G$ and components $\phi_{(u,\j)}=z_1^{j_1}\cdots z_r^{j_r} x_u$ is an eigenvector of the lift graph $G^{\alpha}$ corresponding to the eigenvalue $\lambda$.
		Moreover, all the eigenvalues (including multiplicities) of $G^{\alpha}$ are obtained:
		$$
		\spec G^{\alpha} = \bigcup_{z_1\in R(n_1),\ldots,z_r\in R(n_r)}\spec(\B(\z)).
		$$
	\end{theorem}
	For more information on lift graphs and digraphs, see 
	Dalf\'o, Fiol, Pavl\'ikov\'a, and  \v{S}ir\'an \cite{dfps23}, 
	and Dalf\'o, Fiol, and \v{S}ir\'a\v{n} \cite{dfs19}.

	\section{Johnson graphs as lifts of voltage graphs}
	\label{sec:Johnson}
	
	In this section, we show that the Johnson graphs with some given parameters can be obtained as lifts of base graphs on cyclic groups.
	
	\subsection{Johnson graphs $J(n,2)$ as lifts}
	
	We begin our study by considering the $2$-token graphs of the complete graph $K_n$, where $n$ is an odd number. These graphs correspond to the Johnson graphs $J(n,2)$ or line graphs $L(K_n)$.
	
	\begin{theorem}
		\label{prop:J(n,2)}
		Let $n$ be an odd number, $n=2\nu+1$. Then, the Johnson graph $J(n,2)$ is a lift of a base graph $G$ with $\nu$ vertices with voltages on the group $\Z_n$. 
	\end{theorem}
	
	\begin{proof}
		As representatives of the vertices of $J(n,2)$, that is, vertices of the base graph $G$, we choose 
		$\{0,i\}$ with $i=1,2,\ldots \nu$.
		To determine the voltages, we consider the vertices adjacent to $\{0,i\}$, written in terms of the representatives.
		This gives:
		\begin{align*}
			\{0,i\}  & \sim  \{0,1\},  \{0,2\},\ldots,  \{0,i-1\},\{0,i+1\}, \{0,i+2\}, \ldots,  \{0,\nu\}\\
			& \{0,\nu+1\}=\{0,\nu\}-\nu, \\
			& \{0,\nu+2\}=\{0,\nu-1\}-(\nu-1),\\
			& \vdots \\
			&\{0,n-1\}=\{0,1\}-1,\\
			& \{i,1\}=\{0,i-1\}+1,  \{1,2\}=\{0,i-2\}+2, \ldots,  \{0,i-1\}=\{0,1\}+(i-1),\\
			& \{i,i+1\}=\{0,1\}+i,\ldots,  \{i,i+\nu\}=\{0,\nu\}+i, \\
			& \{i,i+\nu+1\}=\{0,\nu\}-(\nu-i),\\
			& \{i,i+\nu+2\}=\{0,\nu-1\}-(\nu-i-1),\\
			& \vdots \\
			&\{i,n-1\}=\{0,i+1\}-1.
		\end{align*}
		From this, we have that the entries of the row $\{0,i\}$ of the polynomial matrix $\B(z)$ are given by the sums of the elements of the columns of the following array: 
		
		$$
		\begin{array}{cccccccccc}
			(1) & (2) & \cdots & (i-2) & (i-1) & \mbox{\boldmath $(i)$} & (i+1) & \dots & (\nu-1) & (\nu)\\
			1 & 1 & \cdots & 1 & 1 & \mbox{\boldmath $0$} & 1 & \cdots & 1 & 1\\
			\frac{1}{z} & \frac{1}{z^2} & \cdots &\frac{1}{z^{i-2}} &\frac{1}{z^{i-1}} & \mbox{\boldmath $\frac{1}{z^i}$} & \frac{1}{z^{i+1}} &\cdots &\frac{1}{z^{\nu-1}} &\frac{1}{z^{\nu}} \\[0.1cm]
			z^{i-1} & z^{i-2}  & \cdots & z^2 & z & \mbox{\boldmath $0$} &\frac{1}{z} &\cdots &\frac{1}{z^{\nu-i-1}} &\frac{1}{z^{\nu-i}} \\[0.1cm]
			z^i &  z^i & \cdots & z^i & z^i & \mbox{\boldmath $z^i$} &  z^i & \cdots & z^i & z^i 
		\end{array}
		$$
		
		Notice that the elements of the diagonal of $\B(z)$
		are $z^i+\frac{1}{z^i}$ for $i=1,\ldots,\nu$.
		
		Now, let us see that the constructed base graph $G=(V,E)$ with voltages $\alpha:E\rightarrow \Z_n$ forms a lift $G^{\alpha}$ that is isomorphic to  $F_2(K_n)$, under the isomorphism $(\{0,i\},g)\mapsto \{g,g+i\}$.
		To this end, assume that, in $G$, the vertex $\{0,i\}$ is adjacent to  $\{0,j\}$ through an arc with voltage $t\in \Z_n$. This means that, for some $r$, one of the following equalities hold:
		\begin{itemize}
			\item[$(i)$]
			$\{r,i\}=\{0,j\}+t=\{t,t+j\}$;
			\item[$(ii)$]
			$\{0,i+r\}=\{t,t+j\}$.
		\end{itemize}
		Then, in $G^{\alpha}$, the vertex $(\{0,i\},g)\mapsto \{g,g+i\}$ is adjacent to the vertex $(\{0,j\},g+t)\mapsto \{g+t,g+t+j\}$. To prove that, in $F_2(K_n)$, a vertex adjacent to $\{g,g+i\}$
		is $\{g+t,g+t+j\}$, we distinguish the possible cases: 
		If in $(i)$, $r=t$ and $i=t+j$, then $\{g+t,g+t+j\}=\{g+t,g+i\}$;
		on the contrary, if $r=t+j$ and $i=t$, then $\{g+i,g+t+j\}=\{g+t,g+i\}$.
		Similarly if, in $(ii)$, $t=0$ and $i+r=t+j$, then 
		$\{g+t,g+t+j\}=\{g,g+i+r\}$; and, finally, if $t+j=0$ and $i+r=t$, then
		$\{g+t,g+t+j\}=\{g+i+r,g\}$. Thus, in all cases $\{g+t,g+i\}$ or
		$\{g+i+r,g\}$ is a vertex adjacent to $\{g,g+i\}$, as claimed.
	\end{proof}
	
	For instance, for $n=5$ and $7$, the corresponding polynomial matrices are, respectively:
	\begin{equation}
		\B(z)=
		\left(
		\begin{array}{cc}
			z+\frac{1}{z} & 1+z+\frac{1}{z}+\frac{1}{z^2}\\[.1cm]
			1+z+\frac{1}{z}+z^2 & z^2+\frac{1}{z^2}
		\end{array}
		\right),
		\label{B5}
	\end{equation}
	and
	\begin{equation}
		\B(z) = \left(
		\begin{array}{ccc}
			z + \frac{1}{z} & 1 + z + \frac{1}{z} + \frac{1}{z^2} & 1 + z + \frac{1}{z^2} + \frac{1}{z^3} \\[.1cm]
			1 + z + \frac{1}{z} + z^2 & z^2 + \frac{1}{z^2} & 1 + \frac{1}{z} + z^2 + \frac{1}{z^3} \\[.1cm]
			1 + \frac{1}{z} + z^2 + z^3 & 1 + z + \frac{1}{z^2} + z^3 & z^3 + \frac{1}{z^3} \end{array}
		\right).
		\label{B7}
	\end{equation}
	
\begin{figure}[!ht]
		\begin{center}
			\includegraphics[width=12cm]{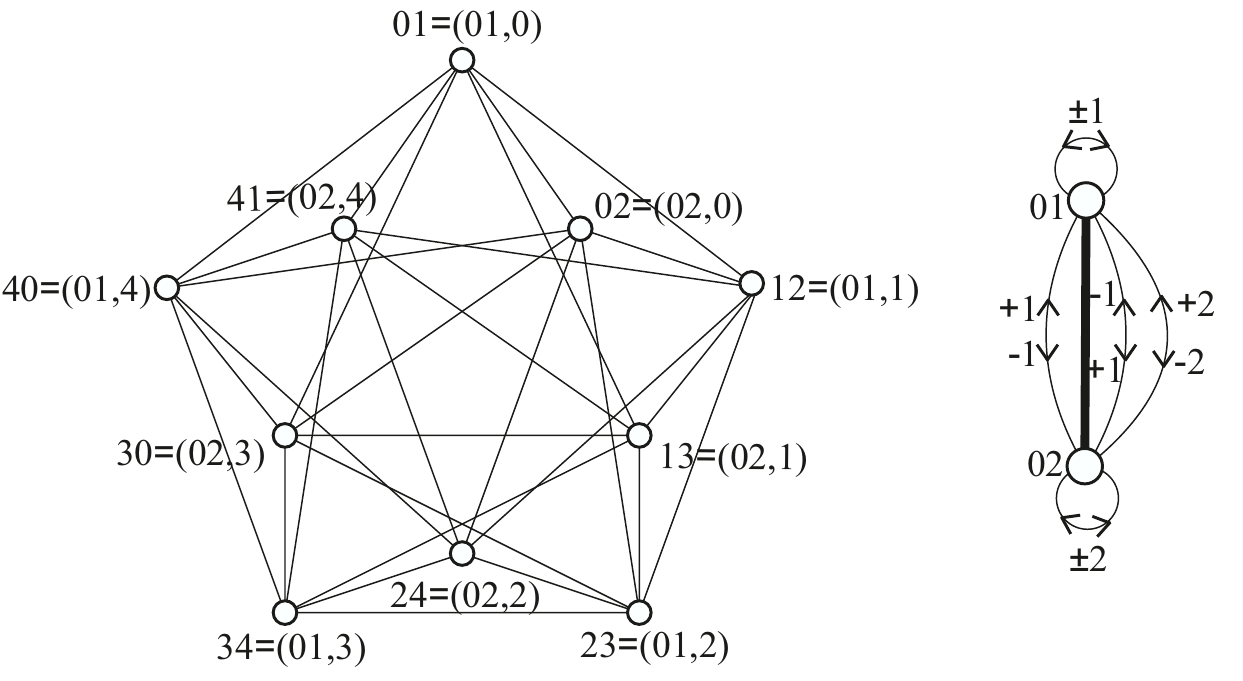}
			\caption{The Johnson graph $J(5,2)$ and its base graph on $\Z_5$. The thick line represents the edge with voltage $0$.}
			\label{fig2}
		\end{center}
	\end{figure}

\begin{figure}[!ht]
		\begin{center}
			\includegraphics[width=5cm]{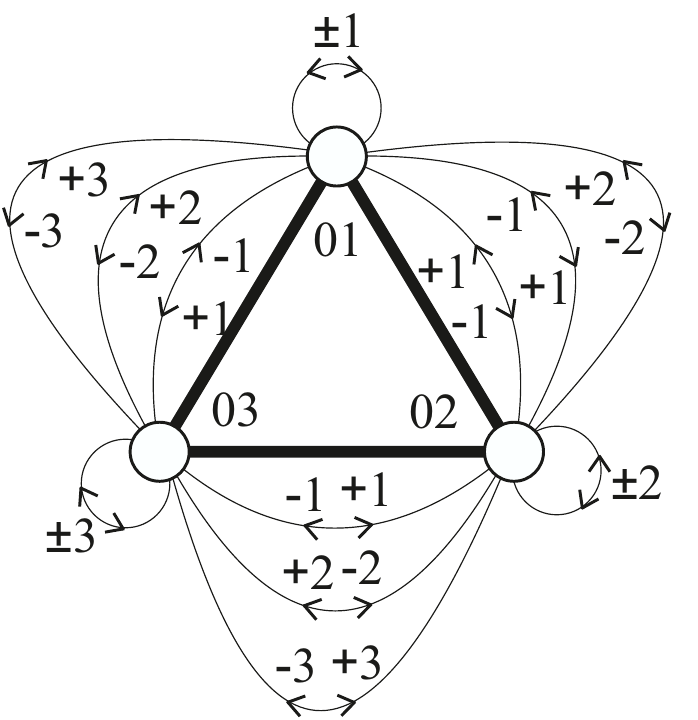}
			\caption{The base graph of $J(7,2)$. The thick lines represent the edges with voltage $0$.}
			\label{fig3}
		\end{center}
	\end{figure}

	Then, for $n=5$ and $n=7$, the Johnson graphs $J(n,2)\cong F_2(K_n)$ are a lift on the cyclic group $\Z_n$ of the base graphs shown in Figures \ref{fig2} and \ref{fig3}.
	In the case of $J(5,2)$, the vertices are labeled as a 2-token of $K_5$, with vertices $0,1,\ldots,4$, and as the lift of the corresponding base graph with voltages on $\Z_5$ (thick lines correspond to voltage $0$).
	
	Moreover, from the matrices \eqref{B5} and \eqref{B7}, we obtain the spectra of $J(n,2)$ for $n=5$ and $n=7$ shown in Tables \ref{table1} and \ref{table2}.
	
	\begin{table}[!ht]
		\begin{center}
			\begin{tabular}{|c|cc| }
				\hline
				$\omega=e^{i\frac{2\pi}{5}}$, $z=\omega^r$ & $\lambda_{r,1}$  & $\lambda_{r,2}$  \\
				\hline\hline
				$\spec(\B(\omega^0))$ & $6$ & $-2$   \\
				\hline
				$\spec(\B(\omega^1))=\spec(\B(\omega^4))$ & $1$ & $-2$   \\
				\hline
				$\spec(\B(\omega^2))=\spec(\B(\omega^3))$ & $1$  &  $-2$ \\
				\hline
			\end{tabular}
		\end{center}
		\caption{All the eigenvalues of the matrices $\B(\omega^r)$, which yield the eigenvalues of the Johnson graph $J(5,2)=F_2(K_5)$.}
		\label{table1}
	\end{table}
	
	\begin{table}[!ht]
		\begin{center}
			\begin{tabular}{|c|ccc| }
				\hline
				$\omega=e^{i\frac{2\pi}{7}}$, $z=\omega^r$ & $\lambda_{r,1}$  & $\lambda_{r,2}$ & $\lambda_{r,3}$ \\
				\hline\hline
				$\spec(\B(\omega^0))$ & $10$ & $-2$ &  $-2$ \\
				\hline
				$\spec(\B(\omega^1))=\spec(\B(\omega^6))$ & $3$ & $-2$  & $-2$ \\
				\hline
				$\spec(\B(\omega^2))=\spec(\B(\omega^5))$ & $3$  &  $-2$  & $-2$\\
				\hline
				$\spec(\B(\omega^3))=\spec(\B(\omega^4))$ & $3$  &  $-$2  & $-2$\\
				\hline
			\end{tabular}
		\end{center}
		\caption{All the eigenvalues of the matrices $\B(\omega^r)$, which yield the eigenvalues of the Johnson graph $J(7,2)=F_2(K_7)$.}
		\label{table2}
	\end{table}
	
	\subsection{Line graphs of circulant graphs as lifts}
	
	In the following result, we show that, in fact, the polynomial matrix $\B(z)$ (associated with  $J(n,2)$ with odd $n$) contains information about the polynomial matrices of the line graphs of circulant graphs, which are Cayley graphs of the cyclic group $\Z_m$ with generators within the interval $[-(m-1), m-1]$.
	
	\begin{theorem}
		Let $n$ be an odd number, $n=2\nu+1$. Let $\B(z)$ be the polynomial matrix of the Johnson graph $J(n,2)\cong F_2(K_n)\cong L(K_n)$, with rows and columns indexed by $\{0,1\},\{0,2\},\ldots,
		\{0,\nu\}$.
		Let $\B'(z)$ be the principal submatrix of $\B(z)$ with rows and columns indexed by $\{0,a_1\},\{0,a_2\},\ldots,\{0,a_s\}$
		with $1\le a_1 < a_2 < \cdots < a_s\le \nu$. Consider $\B'(z)$ as the polynomial matrix of a base graph on $\Z_m$ with $m\ge n$. Then, the lift graph defined by $\B'(z)$ is isomorphic to the line graph $L(H)$, where $H$ is the circulant graph
		$\Cay(\Z_m; \pm a_1,\ldots,\pm a_s)$.
	\end{theorem} 
	\begin{proof}
		Reasoning as at the beginning of this section, notice that the representatives of the vertices of $L(H)$, as vertices of the base graph, can be chosen as $0a_1,0a_2,\ldots,0a_s$. Then, the reasoning is completely similar if we `replace' the elements $1,2,\ldots$ by $a_1,a_2,\ldots$  Indeed, the vertices of $L(H)$ adjacent to $\{0,a_i\}$, written in terms of the representatives, are:
		\begin{align*}
			\{0,a_i\}  & \sim  \{0,a_1\},  \{0,a_2\},\ldots,  \{0,a_{i-1}\},\{0,a_{i+1}\}, \ldots,  \{0,a_{s}\},\\
			& \{0,-a_s\}=\{0,a_s\}-a_s, \ldots,
			\{0,-a_1\}=\{0,a_1\}-a_1,\\
			& \{a_i,a_i-a_1\}=\{0,a_1\}+(a_i-a_1),  \{a_i,a_i-a_2\}=\{0,a_2\}+(a_i-a_2), \ldots \\
			& \{a_i,a_i-a_{i-1}\}=\{0,a_{i-1}\}+(a_i-a_{i-1}),
			\{a_i,a_i-a_{i+1}\}=\{0,a_{i+1}\}+(a_i-a_{i+1}),\ldots\\  &  \{a_i,a_i-a_s\}=\{0,a_{r}\}+(a_i-a_{r}), \\
			& \{a_i,a_i+a_1\}=\{0,a_1\}+a_i,\ldots,
			\{a_i,a_i+a_s\}=\{0,a_s\}+a_i.
		\end{align*}
		Therefore, the entries of the row $\{0,a_i\}$ of the polynomial matrix $\B'(z)$ are given by the sums of the elements of the columns of the following array: 
		$$
		\begin{array}{ccccccccc}
			(a_1) & (a_2) & \cdots & (a_{i-1}) & \mbox{\boldmath $(a_i)$} & (a_{i+1}) & \dots & (a_{s-1}) & (a_s)\\
			1 & 1 & \cdots &  1 & \mbox{\boldmath $0$} & 1 & \cdots & 1 & 1\\
			\frac{1}{z^{a_1}} & \frac{1}{z^{a_2}} & \cdots &\frac{1}{z^{a_{i-1}}} & \mbox{\boldmath $\frac{1}{z^{a_i}}$} & \frac{1}{z^{a_{i+1}}} &\cdots &\frac{1}{z^{a_{s-1}}} &\frac{1}{z^{a_{s}}} \\[.1cm]
			z^{a_i-a_1} & z^{a_i-a_2}  & \cdots &  z^{a_i-a_{i-1}} & \mbox{\boldmath $0$} &\frac{1}{z^{a_{i+1}-a_i}} &\cdots &\frac{1}{z^{a_s-a_i}} &\frac{1}{z^{a_s-a_i}} \\[.1cm]
			z^{a_i} &  z^{a_i} & \cdots  & z^{a_i} & \mbox{\boldmath $z^{a_i}$} &  z^{a_i} & \cdots & z^{a_i} & z^{a_i}\\ 
		\end{array}
		$$
		
		These correspond to the entries of the row indexed by $\{0,a_i\}$ of the principal submatrix $\B'(z)$ of $\B(z)$, as claimed.
	\end{proof}
	
	For instance, if $\B(z)$ is the polynomial matrix of the base graph of Figure \ref{fig3} and $\B'(z)=\B(z)$, the spectra of the  graphs obtained by different values of $m$ are:
	\begin{itemize}
		\item
		$m=7$:
		$\{10, 3^{[6]}, 
		-2^{[14]}\}$ (strongly regular graph $J(7,2)\cong F_2(K_7) \cong L(K_7)$).
		\item
		$m=8$:
		$\{10, 4^{[4]}, 2^{[3]}, 
		-2^{[16]}\}$ (walk regular graph isomorphic to the line graph of $\Cay(\Z_8, \{\pm 1,\pm 2,\pm 3\})$).
		\item
		$m=11$:
		$\{10, 7.148^{[2]}, 4^{[2]}, 2.485^{[2]},
		1.589^{[2]},
		-2^{[18]}\}$ (graph isomorphic to the line graph of $\Cay(\Z_{11}, \{\pm 1,\pm 2,\pm 3\})$).
	\end{itemize}
	Moreover,  the $(2,3)$-submatrix $\B'(z)$ of the same matrix $\B(z)$, that is,
	\begin{equation}
		\B'(z)=\left(
		\begin{array}{cc}
			z^2 + \frac{1}{z^2} & 1 + \frac{1}{z} + z^2 + \frac{1}{z^3} \\[.1cm]
			1 + z + \frac{1}{z^2} + z^3 & z^3 + \frac{1}{z^3} \end{array}
		\right)
		\label{B'12}
	\end{equation}
	corresponds to the base graph of $L(H)$ with $H=\Cay(\Z_m;\pm 2,\pm 3)$. Then if, for example, we take $m=12$, the eigenvalues of $\B'(z)$ with $z=e^{i\frac{r2\pi}{12}}$, for $r=0,\ldots,11$,
	yield
	$$
	\spec(L(H))=\{6^{[1]},3^{[6]},2^{[1]},0^{[2]},-2^{[12]}\}.
	$$
	
	\subsection{Johnson graphs $J(n,k)$ as lifts}
	The following result is a generalization of Theorem \ref{prop:J(n,2)} for $k\ge 2$.
	
	\begin{theorem}
		\label{prop:J(n,k)}
		If $n$ and $k(\le n-k)$ are relatively prime integers, then the Johnson graph $J(n,k)$ is a lift of a base graph $G$ with $\nu={n\choose k}/n$ vertices with voltages on the group $\Z_n$. 
	\end{theorem}
	\begin{proof}
		As before, the proof consists of first identifying the representatives of the orbits (vertices of the base graph), and second defining their `weighted adjacencies' (voltages). Since the procedure to find the latter is the same, we only need to take the case of the former.
		Now, the representatives of the orbits correspond to distinct necklaces with $k$ black beads (representing the $k$ chosen elements of  $\Z_n$) and $n-k$ white beads (corresponding to the remaining elements of $\Z_n$). Thus, two necklaces are equivalent if a given rotation of the other can obtain one. Moreover, if $\gcd(n,k)=1$, all the orbits have $n$ elements, as required.
		These correspond to `aperiodic' necklaces (not consisting of a repeated sequence), also called Lyndon words, see Ruskey and  Sawada \cite{rs99}
		for an efficient algorithm to generate them.
	\end{proof}
	
	For example, to deal with $J(7,3)\cong F_3(K_7)$, and using the simplified notation $ijk$ for $\{i,j,k\}$, we can take the representatives
	$012$, $013$, $014$, $015,$ and $024$. (That is, of each equivalence class, we take the lexicographical smallest sequence).
	Then, the adjacencies of $012$ (the others are similar) in terms of the representatives are:
	\begin{align*}
		012  \sim &\ 013,\  014,\  015,\ 016=012-1,\\
		&\,\, 123=012+1,\ 124=013+1,\ 125=014+1,\ 126=015+1,\\
		&\,\, 023=015+2,\  024,\ 025=024-2,\ 026=013-1.
	\end{align*}
	The $5\times 5$ polynomial matrix, indexed by the above representatives, is
	\begin{equation}
		\B(z)=\left(
		\begin{array}{ccccc}
			z+\frac{1}{z} & 1+z+\frac{1}{z} & 1+z & 1+z+z^2 & 1+\frac{1}{z^2}\\[.1cm]
			1+z+\frac{1}{z} & 0 & 1+\frac{1}{z}+z^3 & 1+z^2+z^3 & z+\frac{1}{z}+z^3\\[.1cm]
			1+\frac{1}{z} & 1+z+\frac{1}{z^3} & z^3+\frac{1}{z^3} & 1+\frac{1}{z}+z^3 & 1+\frac{1}{z^3}\\[.1cm]
			1+\frac{1}{z}+\frac{1}{z^2} & 1+\frac{1}{z^2}+\frac{1}{z^3} & 1+z+\frac{1}{z^3} & 0 & z+\frac{1}{z^2}+z^3\\[.1cm]
			1+z^2 & z+\frac{1}{z}+\frac{1}{z^3} & 1+z^3 & \frac{1}{z}+z^2+\frac{1}{z^3} & z^2+\frac{1}{z^2}
		\end{array}
		\right),
	\end{equation}
	giving, for $z=e^{i\frac{r2\pi}{7}}$ and $r=0,1,\ldots,6$,
	the eigenvalues of $J(7,3)=F_3(K_7)$ shown in Table \ref{table7}.
	
	\begin{table}[!ht]
		\begin{center}
			\begin{tabular}{|c|ccc| }
				\hline
				$\omega=e^{i\frac{2\pi}{7}}$, $z=\omega^r$ & $\lambda_{r,1}$  & $\lambda_{r,2}$ & $\lambda_{r,3}$ \\
				\hline\hline
				$\spec(\B(\omega^0))$ & 12 & $0^{[2]}$ &  $-3^{[2]}$ \\
				\hline
				$\spec(\B(\omega^1))=\spec(\B(\omega^6))$ & 5 & $0^{[2]}$  & $-3^{[2]}$ \\
				\hline
				$\spec(\B(\omega^2))=\spec(\B(\omega^5))$ &  5 & $0^{[2]}$  & $-3^{[2]}$ \\
				\hline
				$\spec(\B(\omega^3))=\spec(\B(\omega^4))$ &  5 & $0^{[2]}$  & $-3^{[2]}$ \\
				\hline
			\end{tabular}
		\end{center}
		\caption{All the eigenvalues of matrices $\B(\omega^r)$, which yield the eigenvalues of the Johnson graph $J(7,3)=F_3(K_7)$.}
		\label{table7}
	\end{table}

	\section{Token graphs of Cayley graphs}
	\label{sec:cay}
	
	Let $\G$ be a group of $n$ elements. For some $k$,
	with $1\le k\le n$, consider the set $\G^{(k)}$ whose elements are the ${n\choose k}$ $k$-subsets $\vecgamma=\{g_1,\ldots,g_k\}$ with $g_i\in \Gamma$. A {\em $k$-set decomposition} ${\cal A}^{(k)}$ of $\G$ is a partition of $\G^{(k)}$ into $r={n\choose k}/n$ classes $A_1,\ldots, A_r$ with $n$ elements each, such that each class $A_i$ is a translation of some  $\vecalpha_i\in A_i$ (with $\alpha_i$ being a `representative' of the class) by $\G$. That is, $A_i=\vecalpha_i\G$. In other words, the classes are the orbits (with the same order $n$) obtained when $\G$ acts on the representatives $\vecalpha_1,\ldots, \vecalpha_r$.
	Thus, a necessary condition for a group $\G$ to have a $k$-set decomposition is that $n$ must divide ${n\choose k}$. For instance, this is the case when dealing with the Abelian group $\Z_{2\mu+1}\times \Z_{2\nu+1}$ and $k=2$ since, then, with $n=(2\mu+1)(2\nu+1)$, we have that the number of classes is $r={n\choose 2}/n=\frac{n-1}{2}=\mu\nu+ \mu+\nu$. In this case, a possible choice of the $r$ 2-sets representatives, of the form $\{(0,0),(i,j)\}=\{00,ij\}$, is as follows:
	$$
	\begin{array}{cccccccc}
		\{00, 0n\} & \{00, 1n\} &  \{00, 2n\} &  \ldots & \{00,mn\} & \{00,m+1,n\} & \ldots & \{00,2m+1,n\}\\
		\vdots &  \vdots  &  \vdots &   &\vdots &\vdots &   &\vdots \\
		{ \bf\{00,01\}} & {\bf\{00,11\}} & {\bf\{00, 21\}} &  \ldots & \{00,m1\} & \{00,m+11\} & \ldots & (00,2m+1,1\}\\
		& {\bf \{00, 10\}} & \{00, 20\} & \ldots & \{00, m0\}&  &  & 
	\end{array}
	$$
	Thus, in particular, the group
	$\G=\Z_3\times \Z_3$ has the following 2-set decomposition with four orbits having representatives
	$\vecalpha=\{00,01\}$, $\vecbeta=\{00,10\}$, $\vecgamma=\{00,11\}$, and $\vecdelta=\{00,21\}$ (above, in boldface). The elements of the same orbit are indicated by the same Greek letter as their representative. For instance, the orbit of $\vecalpha=(00,10)$ has elements $(00,20)$, $(10,20)$, $(01,11)$, \ldots (see Table \ref{table4}).
	
	\begin{table}[!ht]
		\begin{center}
			\begin{tabular}{c|ccccccccc|}
				& 00 & 10 & 20 & 01 & 11 & 21 & 02 & 12 & 22\\
				\hline
				00 &  & \mbox{{\boldmath $\alpha$}} & $\alpha$ & $\vecbeta$ & $\vecgamma$ & $\vecdelta$ & $\beta$ & $\delta$ & $\gamma$  \\ 
				10 &  &  & $\alpha$ & $\delta$ & $\beta$ & $\gamma$  & $\gamma$ & $\beta$ & $\delta$  \\ 
				20 &  &  &  & $\gamma$ & $\delta$ & $\beta$ & $\delta$ & $\gamma$ & $\beta$  \\ 
				01 &  &  &  &  & $\alpha$ & $\alpha$ & $\beta$ & $\gamma$ & $\delta$  \\ 
				11 &  &  &  &  &  & $\alpha$ & $\delta$ & $\beta$ & $\gamma$  \\ 
				21 &  &  &  &  &  &  & $\gamma$ & $\delta$ & $\beta$  \\ 
				02 &  &  &  &  &  &  &  & $\alpha$ & $\alpha$  \\ 
				12 &  &  &  &  &  &  &  &  & $\alpha$  \\ 
				22 &  &  &  &  &  &  &  &  &   \\ 
				\hline 
			\end{tabular}
		\end{center}
		\caption{The elements of the same orbit are indicated by the same Greek letter as their representative in the group $\G=\Z_3\times \Z_3$. The boldface indicates the representatives.}
		\label{table4}
	\end{table}
	
	\begin{theorem}
		Let $\G$ be a group with order $n$.
		For some generating set $S=\{\pm a_1,\ldots,\pm a_s\}$ of $\G$, consider the Cayley graph $G=\Cay(\G, S)$, and its $k$-token graph $F_k(G)$ for some $k$. If $\G$ has a $k$-set decomposition, then $F_k(G)$ can be obtained as the lift of a base graph $H$ with $r={n\choose k}/n$ vertices and voltages on the group $\G$.
	\end{theorem}
	\begin{proof}
		As in the case of cyclic groups, the vertices of $H$ are labeled with the $r$ representatives of the $k$-set decomposition.
		Now to fix the voltages, we reason as before: If vertex $\vecalpha=\{\alpha_1,\alpha_2,\ldots,\alpha_k\}$ is adjacent to, say, $\vecalpha'=\{\alpha_1a_i,\alpha_2,\ldots,\alpha_k\}$ for some $a_i\in S$, there must be some $g\in \Gamma$ and $\vecbeta\in {\cal A}^{(k)}$ such that $\vecalpha'=\vecbeta g$. Then, the arc in $H$ from $\vecalpha$ to $\vecbeta$
		has voltage $g$.
	\end{proof}
	
\begin{figure}[!ht]
		\begin{center}
			\includegraphics[width=6cm]{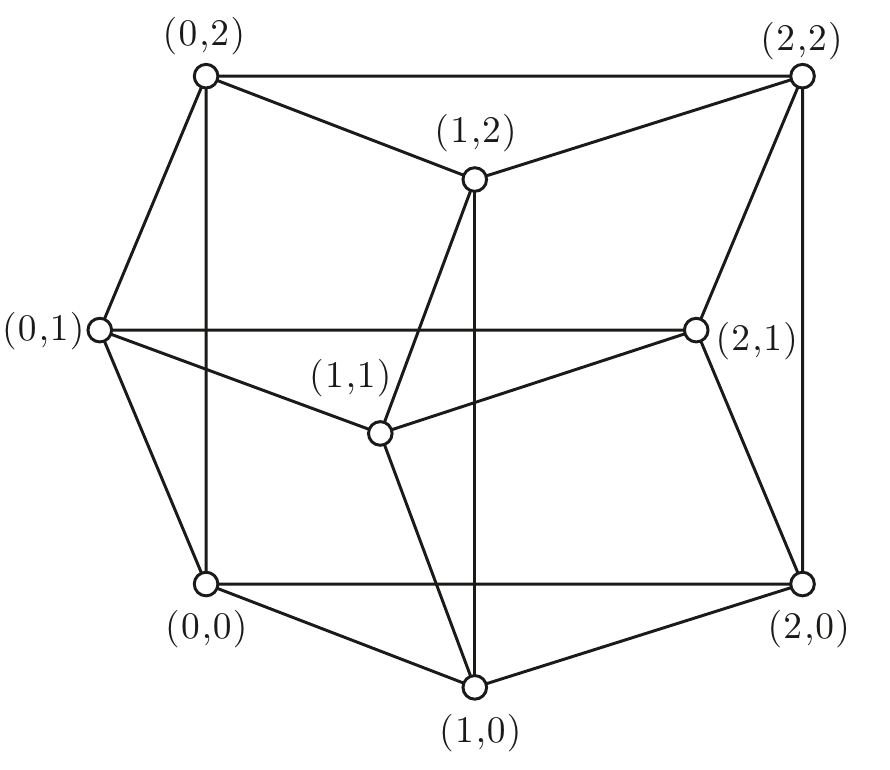}
			\caption{The Cayley graph $\Cay(\G,S)$ with $\G=\Z_3\times \Z_3$ and $S=\{\pm (1,0),\pm (0,1)\}$.}
			\label{fig4}
		\end{center}
	\end{figure}

	Let us show an example with the 2-token graph of $\Cay(\G,S)$ with $\G=\Z_3\times \Z_3$ and $S=\{\pm (1,0),\pm (0,1)\}$, see Figure \ref{fig4}.
	
	\begin{align*}
		\vecalpha=\{00,10\}  \sim &\  \{00,20\}=\vecalpha+20, \{00,11\}=\vecgamma, \{00,12\}=\vecdelta+12,\\
		& \ \{20,10\}=\vecalpha+01, \{01,10\}=\vecdelta+10, \{02,10\}=\vecgamma+02.\\
		\vecbeta=\{00,01\}  \sim &\  \{00,11\}=\vecgamma, \{00,21\}=\vecdelta, \{00,02\}=\vecbeta+02,\\
		& \ \{10,01\}=\vecdelta+10, \{20,01\}=\vecgamma+20, \{02,01\}=\vecbeta+01.\\                  
		\vecgamma=\{00,11\}  \sim &\  \{00,21\}=\vecdelta, \{00,01\}=\vecbeta, \{00,12\}=\vecdelta+12, \{00,10\}=\vecalpha,\\
		& \ \{10,11\}=\vecbeta+10, \{20,11\}=\vecdelta+20, \{01,11\}=\vecalpha+01, (02,11)=\vecdelta+11.
	\end{align*}
	\begin{align*}
		\vecdelta=\{00,21\}  \sim &\  \{00,01\}=\vecbeta, \{00,11\}=\vecgamma, \{00,22\}=\vecgamma+22, \{00,20\}=\vecalpha+20,\\
		& \ \{10,21\}=\vecgamma+10, \{20,21\}=\vecbeta+20, \{01,21\}=\vecalpha+21, \{02,21\}=\vecgamma+21.
	\end{align*}
	
	Then, the (two variables) polynomial matrix of the 2-token of the `toroidal mesh' $G=\Cay(\G,\{\pm 01,\pm 10\})$ is
	$$
	\B(y,z)=\left(
	\begin{array}{cccc}
		y+y^2  & 0 & 1+z^2 & y+yz^2\\
		0 & z+z^2 & 1+y^2 & 1+y\\
		1+z & 1+y & 0 & 1+y^2+yz+yz^2\\
		y^2+y^2z & 1+y^2 & 1+y+y^2z+y^2z^2 & 0
	\end{array}
	\right).
	$$
	With the 9 possible pairs of values $(y,z)$, for $y=e^{i\frac{r2\pi}{3}}$, $z=e^{i\frac{s2\pi}{3}}$, and
	$r,s=0,1,2$, the obtained (approximated) eigenvalues are shown in Table \ref{tab:abelian}.
	
	\begin{table}
		\begin{center}
			\begin{tabular}{|c|ccc|}
				\hline
				$r \backslash s$  & 0 & 1 & 2 \\
				\hline
				\hline
				0 & 7.12, 1.34,$-1$,$-3.13$ & 3.78, 1.54, $-1$, $-3.13$ & 3.78, 1.54, $-1$, $-3.13$ \\ 
				1 & 3.78, 1.54, $-1$, $-3.1$3 & 2, 0.56, $-1$, $-3.56$ &  2, 0.56, $-1$, $-3.5$6\\ 
				2 & 3.78, 1.54, $-$1, $-3.13$ &  2, 0.56, $-1$, $-3.56$ &  2, 0.56, $-1$, $-3.56$ \\ 
				\hline 
			\end{tabular}
		\end{center}
		\caption{The eigenvalues of the $2$-token graph of $\Cay(\G, S)$ with $\G=\Z_3\times \Z_3$ and $S=\{\pm 10, \pm 01\}$.}
		\label{tab:abelian}
	\end{table}
	
	In the general case, when the group $\G$ is not necessarily Abelian, we proceed as follows. If $G=(V,E)$ is the base graph with 
	voltage assignment $\alpha:E \rightarrow \G$, its `base matrix' $\B(G)$ has entries $\B(G)_{uv}=\alpha (a_1)+\cdots +\alpha (a_r)\in \CC[\G]$, where $a_1,\ldots, a_r$ are the arcs from $u$ to $v$, and $\B(G)=0$ if $(u,v)\not\in E$. 
	Then, if $\rho\in \Irep(\G)$ is a unitary irreducible representation of $\G$, with dimension $d_{\rho}=\dim(\rho)$, the $d_{\rho}n\times d_{\rho}n$ matrix $\rho(\B(G))$ (where $\rho$ acts on the elements of $\G$) plays the same role as the polynomial matrix $\B(\z)$ with a given value of $\z$.
	To find the spectrum of the lift $G^{\alpha}$, we apply the following result, which can be seen as a generalization of Theorem \ref{th:sp-lifts}.
	
	\begin{theorem}[\cite{dfs19}]
		\label{theo-sp}
		Let $G=(V,E)$ be a base (di)graph on $n$ vertices, with a voltage assignment $\alpha$ in a group $\G$.
		For every irreducible representation $\rho\in \Irep (\G)$, let $\rho(\B)$ be the complex matrix whose entries are given by $\rho(\B(G)_{u,v})$.
		Then,
		$$
		\spec G^{\alpha} = \bigcup_{\rho\in \Irep(\G)}d_\rho\cdot\spec(\rho(\B(G))).
		$$
	\end{theorem}
	
	\section{Further results}
	\label{sec:further}
	
	The results of the previous sections can be easily extended in the following directions.
	
	\subsection{Computing eigenspaces}
	
	The results of Theorems \ref{th:sp-lifts} and \ref{theo-sp} are based on the fact that every eigenvector of the base graph $G$ yields an eigenvector of the lift $G^{\alpha}$. Thus, although not explained in the examples, the eigenspaces of our constructions can also be computed.
	
	\subsection{Universal matrix}
	
	As it was shown by the authors in \cite{dfps23}, from the base graph $G$ with voltages, we can associate a polynomial-type matrix representing the so-called universal adjacency matrix of the lift $G^{\alpha}$. 
	The universal adjacency matrix $\U$ of a graph, with adjacency matrix $\A$, is a linear combination, with real
	coefficients, of $\A$, the diagonal matrix $\D$ of vertex degrees, the identity matrix $\I$, and the all-1 matrix $\J$; that is, $\U=c_1 \A+c_2 \D+c_3 \I+ c_4 \J$, with $c_i\in {\mathbb R}$ and $c_1\neq 0$. Thus, as particular cases, $\U$ may be the adjacency matrix, the Laplacian, the signless Laplacian, and the Seidel matrix (see, for instance, Haemers and Omidi \cite{ho11}).
	Then, our method also gives the eigenvalues and eigenspaces of such matrices.
	
	\subsection{Dealing with digraphs}
	
	As commented in the preliminary section, a graph is a type of digraph in which two opposite arcs constitute each edge.
	In fact, this paper's results and argumentation can be generalized naturally to digraphs.
	For instance, the vertices of the $k$-token digraph $F_k(G)$ of a digraph $G=(V,A)$, on $n$ vertices, can be represented by $k$ indistinguishable tokens placed in different vertices of $G$. One vertex of $F_k(G)$ is adjacent to another one (forming an arc) if the latter can be obtained from the former by moving one token along an arc in $A$.
	
\begin{figure}[!ht]
		\begin{center}
			\includegraphics[width=8cm]{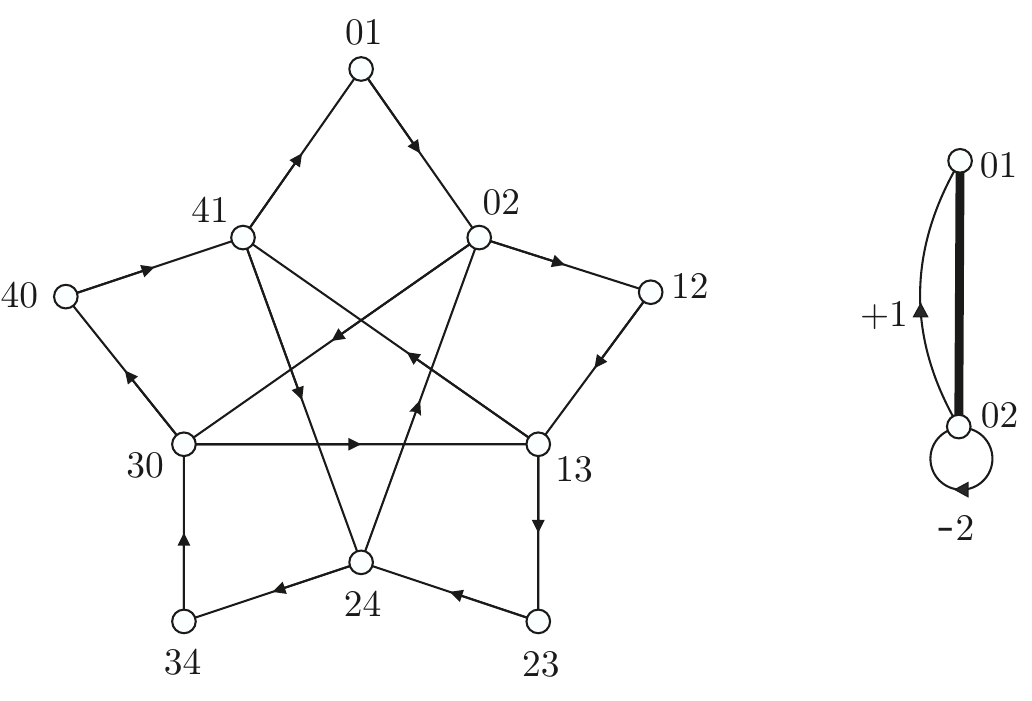}
			\caption{The $2$-token digraph of the directed cycle $C_5$ (with arcs $0\rightarrow 1\rightarrow 2\rightarrow 3\rightarrow 4\rightarrow 0$) and its base digraph. The thick line represents the edge with voltage $0$.}
			\label{fig5}
		\end{center}
	\end{figure}

	Of course, when dealing with digraphs, the eigenvalues and eigenvectors are not necessarily real, but this does not affect the procedure for obtaining them.
	For example, the $2$-token  digraph of the directed cycle $C_5$, shown in Figure 
	\ref{fig5}, can be constructed as a lift of the base digraph on $\Z_5$, with 
	vertices $01$ and $02$ and three arcs: $(01,02)$ with voltage 0; $(02,01)$ with voltage $+1$; and loop $(02,02)$ with voltage $-2$. Then, the polynomial matrix is
	$$
	\B(z)=
	\left(
	\begin{array}{cc}
		0 & 1 \\
		z & \frac{1}{z^2}
	\end{array}
	\right).
	$$
	
	By taking $z=e^{i\frac{r2\pi}{5}}$, for $r=0,\ldots,4$, the eigenvalues of $\B(z)$
	yield all the eigenvalues of $F_2(C_5)$. Namely,
	$$
	1.618,\ 0.5\pm 1.540\, i,\ 0.5\pm 0.363\, i,\ -0.618,\ -0.191\pm 0.588\, i, \ -1.309\pm 0.951\, i.
	$$
	
	\section*{Statements and Declarations}
	
	The authors have no competing interests.
	


\begin{thebibliography}{00}
		
		\bibitem{abel91}
		Y. Alavi, M. Behzad, P. Erd\H{o}s, and D. R. Lick,
		Double vertex graphs
		{\em J. Combin. Inform. System Sci.} {\bf 16} (1991), no. 1, 37--50.
		
		\bibitem{agrr07}
		K. Audenaert, C. Godsil, G. Royle, and T. Rudolph,
		Symmetric squares of graphs,
		\emph{J. Combin. Theory B} \textbf{97} (2007) 74--90.
		
		\bibitem{bbmprs97}
		E. T. Baskoro, L. Brankovi\'c, M. Miller, J. Plesn\'{\i}k, J. Ryan and J. \v{S}ir\'a\v{n}, Large digraphs with small diameter: A voltage assignment approach, 
		{\em J. Combin. Math. Combin. Comput.} {\bf 24} (1997) 161--176.
		
		\bibitem{biggs}
		N. Biggs, 
		{\it Algebraic Graph Theory}, 
		Cambridge Univ. Press, Cambridge, 2nd ed., 1993.
		
		\bibitem{cgss76}
		P. J. Cameron, J. M. Goethals, J. J. Seidel, and E. E. Shult, 
		Line graphs, root systems, and elliptic geometry, 
		{\em J. Algebra} {\bf 43} (1976) 305--327.
		
		\bibitem{ddffhtz21}
		C. Dalf\'o, F. Duque, R. Fabila-Monroy, M. A. Fiol, C. Huemer, A. L. Trujillo-Negrete, F. J. Zaragoza Mart\'inez,
		On the Laplacian spectra of token graphs,
		{\em Linear Algebra Appl.} {\bf 625} (2021) 322--348.
		
		
		\bibitem{dfmrs17}
		C. Dalf\'o, M. A. Fiol, M. Miller, J. Ryan, and J. \v{S}ir\'a\v{n},
		An algebraic approach to lifts of digraphs,
		{\em Discrete Appl. Math.} {\bf 269} (2019) 68--76.
		
		\bibitem{dfps23}
		C. Dalf\'o, M. A. Fiol, S. Pavl\'ikov\'a, and J. \v{S}ir\'an, 
		On the spectra and eigenspaces of the universal adjacency matrices of arbitrary lifts of graphs, 
		\textit{Linear Multilinear Algebra} {\bf 71} (2023), no. 5, 693--710.
		
		\bibitem{dfs19}
		C. Dalf\'o, M. A. Fiol, and J. \v{S}ir\'a\v{n},
		The spectra of lifted digraphs, 
		{\em J. Algebraic Combin.} \textbf{50} (2019) 419--426.
		
		\bibitem{ffhhuw12}
		R. Fabila-Monroy, D. Flores-Pe\~{n}aloza, C. Huemer, F. Hurtado, J. Urrutia, and D. R. Wood, Token graphs,
		\emph{Graphs Combin.} \textbf{28} (2012), no. 3, 365--380.
		
		\bibitem{flt22}
		R. Fabila-Monroy, J. Leaños, and A. L. Trujillo-Negrete, 
		On the connectivity of token graphs of trees, 
		{\em Discrete Math. Theor. Comput. Sci.} {\bf 24} (2022), no. 1, Paper No. 9, 23 pp.
		
		\bibitem{g93}
		C. D. Godsil, 
		{\it Algebraic Combinatorics}, 
		Chapman and Hall, New York, 1993.
		
		
		
		\bibitem{gt87}
		J. L. Gross and T. W. Tucker, 
		{\em Topological Graph Theory}, 
		Wiley, New York (1987).
		
		\bibitem{ho11}
		W. H. Haemers and G. R.  Omidi,
		Universal adjacency matrices with two eigenvalues,
		{\em Linear Algebra Appl.} {\bf 435} (2011), no. 10, 2520--2529.
		
		\bibitem{ms13}
		M. Miller and J. \v{S}ir\'a\v{n},
		Moore graphs and beyond: A survey of the degree/diameter problem, 
		{\em Electron. J. Combin.} {\bf 20(2)} (2013) \#DS14v2.
		
		
		\bibitem{rdf23}
		M. A. Reyes, C. Dalf\'o, and M. A. Fiol,
		On the spectra and spectral radii of token graphs, 
		\textit{Bol. Soc. Mat. Mex.} (2023) (accepted).
		
		\bibitem{rdfm23}
		M. A. Reyes, C. Dalf\'o,  M. A. Fiol, and A. Messegu\'e,
		A general method to find the spectrum and eigenspaces of the $k$-token of a cycle, and
		2-token through continuous fractions, 
		submitted, 2023.
		
		\bibitem{rs99}
		F. Ruskey and J. Sawada.
		An efficient algorithm for generating necklaces with fixed density, 
		{\em SIAM J. Comput.} {\bf 29} (1999), no. 2, 671--684.
		
		
	\end{thebibliography}
\end{document}